\documentclass[12pt,twoside,reqno]{amsart}

\usepackage{amsmath}
\usepackage{amssymb}
\usepackage{latexsym}
\usepackage[all]{xy}

\newtheorem{theorem}{\bf Theorem}

\newtheorem{corollary}[theorem]{\bf Corollary}

\newtheorem{definition}[theorem]{\bf Definition}
\newtheorem{example}[theorem]{\bf Example}
\newtheorem{remark}[theorem]{\bf Remark}

\def\proof{\noindent {\sc Proof:}\enspace}
\def\endproof{\removelastskip\rightline{$\Box$}\par\bigskip}

\def\conj#1{ \{ #1 \} }

\begin{document}
\setcounter{page}{1}
\title{Combinatorial derived invariants for gentle algebras}

\author{Diana Avella-Alaminos}
\address{Diana Avella-Alaminos\newline Departamento de Matem\'aticas \newline Facultad de Ciencias \newline Universidad Nacional Aut\'onoma de M\'exico.
        \newline Ciudad Universitaria\newline M\'exico D.F. C.P. 04510, M\'exico}
\email{avella@matem.unam.mx}

\author{Christof Geiss}
\address{Christof Geiss \newline Instituto de Matem\'aticas
\newline Universidad Nacional Aut\'onoma de M\'exico.
         \newline Ciudad Universitaria\newline
        M\'exico D.F. C.P. 04510, M\'exico}
        \email{christof@matem.unam.mx}

\thanks{The research was supported by DGAPA and DGEP, UNAM and CONACyT}
\date{July 2007}

\begin{abstract}
We define derived equivalent invariants for gentle algebras, constructed in an easy combinatorial way from the quiver with relations defining these algebras. Our invariants consist of pairs of natural numbers and contain important information about the algebra and the structure of the stable Auslander-Reiten quiver of its repetitive algebra.
As a by-product we obtain that the number of arrows of the quiver of a gentle algebra is invariant under derived equivalence. Finally, our invariants separate the derived equivalence classes of gentle algebras with at most one cycle.
\end{abstract}
\maketitle

\section{Introduction}

Let ${\mathrm k}$ be a field. For a finite-dimensional connected
${\mathrm k}$-algebra $A$, the bounded derived category $D^b(A)$
of the module category of finite-dimensional left $A$-modules,
$A$-$\operatorname{mod}$, contains in some sense the complete
homological information about $A$. So it is natural to classify
such algebras up to derived equivalence. Examples of such
invariants are the Grothendieck groups, Hochschild (co-)homology
and cyclic homology or the neutral component of the group of outer
automorphisms, see \cite{hs01} and \cite{ro00}. However it is in
general not very easy to calculate these invariants.

In this paper we focus on gentle algebras. This class of algebras
is closed under derived equivalence \cite{sz}. We propose for
gentle algebras new invariants $\phi_A:\mathbb{N}^2\to\mathbb{N}$
which can be determined easily if $A$ is given as a quiver with
relations  ${\mathrm k}Q/ \left< \mathcal{P} \right>$.

Roughly speaking $\phi_A$ is obtained as follows: Start with a
maximal directed path in $Q$ which contains no relations. Then
continue in opposite direction as long as possible with zero
relations. Repeat this until the first path appears again, say
after $n$ steps. Then we obtain a pair $(n,m)$ where $m$ is the
number of arrows which appeared in a zero relation. Repeat this
procedure until all maximal paths without a zero relation have
been used; $\phi_A$ counts then how often each pair $(n,m)\in
\mathbb{N}^2$ occurred. See Section \ref{algoritmo} for a precise
description.

In order to show that this is in fact a derived invariant, recall
first that the repetitive algebra $\hat{A}$ is special biserial
and selfinjective. We show that $\phi_A$ describes the action of
the suspension functor $\Omega^{-1}_{\hat{A}}$ for the
triangulated category $\hat{A}$-$\underline{\operatorname{mod}}$
on the Auslander-Reiten components which contain string modules
and Auslander-Reiten triangles, see \cite[1.4]{ha88} of the form
$X\to Y\to \tau^{-1}_{\hat{A}}X\to \Omega^{-1}_{\hat{A}}X$ with
$Y$ indecomposable. Note that such components are of the form
$\mathbb{Z}A_{\infty}$ or $\mathbb{Z}A_{\infty}/\left< \tau^r
\right>$ for some $r \geq 1$.

Since $D^b(A)\cong D^b(B)$ implies that also
$\hat{A}$-$\underline{\operatorname{mod}}\cong\hat{B}$-$\underline{\operatorname{mod}}$
as triangulated categories (even if gl.dim$A=\infty$) we conclude
our main result, see Section \ref{justificacion}.\\
%$\phantom {11111111111111111111111111111111111111111111111111111111111111111111111}$\\
%\vspace{.5 cm}

 {\bfseries Theorem A.}\label{soninv} {\em Let $A$ and $B$ be
gentle algebras. If $A$ and $B$ are derived equivalent then
$\phi_A=\phi_B$.}

\vspace{0.5 cm}

We should point out that the argument to show that $\phi_A$ is in
fact a derived invariant, is quite delicate in case of homogeneous
tubes which come from string modules. For this reason we need to
analyze the outer automorphism group of a gentle algebra. As a
by-product we obtain in Section \ref{automorfismos} the following result:\\

{\bfseries Proposition B.}\label{aristasinv} {\em The number of
cycles and the number of arrows are invariant under derived
equivalence for gentle algebras.}

\vspace{0.5 cm}

  Other combinatorial invariants for gentle algebras were recently
obtained in \cite{bh05}. In Section \ref{clas} we apply our
invariant to the known derived classification of gentle algebras
where $Q$ admits at most one cycle (see \cite{ah81}, \cite{as87},
\cite{dv01} and
\cite{bg04}).\\

{\bfseries Theorem C.}\label{unciclo} {\em Let $A={\mathrm k}Q/
\left< \mathcal{P} \right>$ and $B={\mathrm k}Q'/ \left<
\mathcal{P}' \right>$ be gentle algebras such that $c(Q),c(Q')\leq
1$. Then $A$ and $B$ are derived equivalent if and only if
$\phi_A=\phi_B$. }

\vspace{0.5 cm}

In a forthcoming paper, the first named author will prove that if
$A$ is a gentle algebra with two cycles, $\#\text{ Supp } (\phi_A)
 \in \{ 1,3\}$ and will use these invariants to classify those
kind of algebras in the case where $\#\text{ Supp } (\phi_A) =3$.

 The results
presented in this article form part of the Ph.D. thesis of the
first named author. We should like to thank Corina Saenz and
Bertha Tom\'e for useful conversations.

\section{Preliminaries}
\subsection{Gentle algebras}

A {\em quiver} is a tuple $Q=(Q_0,Q_1,s,e)$, where $Q_0$ is the
set of vertices, $Q_1$ the set of arrows and $s,e:Q_1 \rightarrow
Q_0$ the functions which determine the start resp. end point of an
arrow.  For simplicity we assume always that $Q$ is connected.
When $Q$ is finite, we say {\em $Q$ has $c(Q)$ cycles} if $\#Q_1 =
\#Q_0+c(Q)-1$, in other words if $c(Q)$ is the least number of
arrows that we have to remove from $Q$ in order to obtain a tree.
A {\em path} in $Q$ is a sequence of arrows $C=\alpha_n \dots
\alpha_2  \alpha_1$ such that $s(\alpha_{i+1})=e(\alpha_{i})$ for
$1\leq i\lneq n$, its {\em length} is $n$, $l(C):=n$. For $v\in
Q_0$ define a {\em trivial path} $1_v$ of length zero. We can
extend $s$ and $e$ to paths in the obvious way. For $\alpha \in
Q_1$ define $\alpha ^{-1} $ with $s(\alpha^{-1}):=e(\alpha)$,
$e(\alpha ^{-1}):=s(\alpha)$ and $(\alpha ^{-1})^{-1}:=\alpha$;
for a path $C=\alpha_n \dots \alpha_2 \alpha_1$ define
$C^{-1}:=\alpha_1^{-1} \alpha_{2}^{-1} \dots \alpha_n^{-1}$ and
$1_v^{-1}:=1_v$ for a trivial path. If $C_1$ and $C_2$ are paths
in $Q$ define the composition $C_2C_1$ as the concatenation of
both paths if $s(C_2)=e(C_1)$ or $0$ if $s(C_2)\neq e(C_1)$. For a
field ${\mathrm k}$ let ${\mathrm k}Q$ be the corresponding {\em
path algebra}, with paths of $Q$ as a basis and multiplication
induced by composition. A non zero linear combination of paths of
length at least two, with the same start point and end point is
called a {\em relation in $Q$}. Let $\mathcal{P}$ be a set of
relations in $Q$ and $\left< \mathcal{P} \right>$ the ideal of
${\mathrm k}Q$ generated by $\mathcal{P}$. In general we consider
algebras ${\mathrm k}Q/ \left< \mathcal{P} \right>$ and we
identify a path in $Q$ with its corresponding class in ${\mathrm
k}Q/ \left< \mathcal{P} \right>$.

The algebras we are interested in are defined by quivers with relations which fulfill very particular conditions.

\begin{definition}\label{biserial}{\normalfont
We call ${\mathrm k}Q/ \left< \mathcal{P} \right>$ a {\em special
biserial algebra} if the following three conditions hold:
\begin{enumerate}
\item For each $v \in Q_0$, $\#\{ \alpha \in Q_1 | s(\alpha) = v
\} \leq 2$ and $\#\{ \alpha \in Q_1 | e(\alpha) = v \} \leq 2$.
\item For each $\beta \in Q_1$, $\#\{ \alpha \in Q_1 | s(\beta) =
e(\alpha) \text{ and } \beta \alpha \notin \mathcal{P} \} \leq 1$
and \\$\#\{ \gamma \in Q_1 | s(\gamma) = e(\beta) \text{ and }
\gamma \beta \notin \mathcal{P} \} \leq 1$.
\item For each $\beta \in Q_1$ there is a bound $n(\beta)$ such that any path $\beta_{n(\beta)}\dots \beta_2 \beta_1$ with $\beta_{n(\beta)}=\beta $ contains a subpath in $\mathcal{P}$ and any path $\beta_{n(\beta)}\dots\beta_2 \beta_{1}$ with $\beta_{1}=\beta$ contains a subpath in $\mathcal{P}$.\\

$A$ is {\em gentle} if moreover:\\

\item All relations in $\mathcal{P}$ are monomials of length $2$.
\item For each $\beta \in Q_1$, $\#\{ \alpha \in Q_1 | s(\beta) =
e(\alpha) \text{ and } \beta \alpha \in \mathcal{P} \} \leq 1$ and
\\$\#\{ \gamma \in Q_1 | s(\gamma) = e(\beta) \text{ and } \gamma
\beta \in \mathcal{P} \} \leq 1$.
\end{enumerate}
}\end{definition}

\begin{remark}{\normalfont
A special biserial algebra is finite-dimensional if and only if
$Q$ is finite, see \cite[1]{ri97}. }\end{remark}
\begin{remark}{\normalfont

A special biserial algebra is called a {\em string algebra} if
$\mathcal{P}$ consists only of paths. }
\end{remark}
We will make an slight abuse of notation by talking about gentle
algebras referring to the quiver with relations which define those
algebras.

\subsection{Threads of a gentle algebra}
Let $A$ be a string algebra. A {\em permitted path} of $A$ is a path $C=\alpha_n \dots \alpha_2 \alpha_1$ with no zero relations and it is a {\em non trivial permitted thread} of $A$ if for all $\beta\in Q_1$, neither $C\beta$ nor $\beta C$ is a permitted path. Similarly a {\em forbidden path} of $A$ is a sequence $\Pi=\alpha_n \dots \alpha_2 \alpha_1$ formed by pairwise different arrows in $Q$ with $\alpha_{i+1} \alpha_{i} \in \mathcal{P}$ for all $i \in \{ 1,2,\ldots n-1\}$ and it is a {\em non trivial forbidden thread} if for all $\beta\in Q_1$, neither $\Pi\beta$ nor $\beta \Pi$ is a forbidden path. The existence of non trivial permitted threads is due to point $(3)$ in Definition \ref{biserial} and the existence of non trivial forbidden threads is due to the restriction of considering pairwise different arrows.

We would like every vertex to be involved in exactly two permitted threads and exactly two forbidden threads so we define also trivial threads for some vertices. Let $v \in Q_0$ such that $\#\{ \alpha \in Q_1 | s(\alpha)=v  \} \leq 1$,  $\#\{ \alpha \in Q_1 | e(\alpha)=v  \} \leq 1$ and if $\beta ,\gamma \in Q_1$ are such that $s(\gamma)=v=e(\beta)$ then $\gamma \beta \notin \mathcal{P}$, we consider $1_v$ a {\em trivial  permitted thread} in $v$ and denote it by $h_{v}$. Let $\mathcal{H}_A$ be the set of all permitted threads of $A$, trivial and non trivial. This set describes completely the algebra $A$.
Similarly, for $v \in Q_0$ with $\#\{ \alpha \in Q_1 | s(\alpha)=v  \} \leq 1$,  $\#\{ \alpha \in Q_1 | e(\alpha)=v  \} \leq 1$ and if $\beta ,\gamma \in Q_1$ are such that $s(\gamma)=v=e(\beta)$ then  $\gamma \beta \in \mathcal{P}$, we consider $1_v$ a {\em trivial forbidden thread} in $v$ and denote it by $p_{v}$. Observe that certain paths can be permitted and forbidden threads at the same time.

For a string algebra it is possible to describe the relations in its quiver by using two functions $\sigma ,\varepsilon :Q_1 \to \{ 1,-1\} $ as in \cite[3]{br87}, defined by:
\begin{enumerate}
\item
If $\beta _1 \neq \beta _2$ are arrows with $s(\beta _1)=s(\beta _2)$, then $\sigma (\beta _1)= -\sigma (\beta _2)$.
\item
If $\gamma _1 \neq \gamma _2$ are arrows with $e(\gamma _1)=e(\gamma _2)$, then $\varepsilon (\beta _1)=-\varepsilon (\beta _2)$.
\item
If $\beta$, $\gamma$ are arrows with $s(\gamma )=e(\beta )$ and $\gamma \beta \notin \mathcal{P}$, then $\sigma (\gamma)= -\varepsilon(\beta)$.
\end{enumerate}

We can extend these functions to threads of $A$ as follows. For
$H=\alpha_n \dots \alpha_2 \alpha_1$ a non trivial thread of $A$
define $\sigma (H):=\sigma (\alpha_1)$ and
$\varepsilon(H):=\varepsilon(\alpha_n)$. If there is a trivial
permitted thread $h_{v}$ for some $v \in Q_0$, connexity of $Q$
assures the existence of some $\gamma \in Q_1$ with $s(\gamma)=v$
or $\beta \in Q_1$ with $e(\beta)=v$, in the first case we define
$\sigma (h_v)=-\varepsilon (h_v)=-\sigma (\gamma)$ , for the
second case $\sigma(h_v)=-\varepsilon(h_v)=\varepsilon (\beta)$.
If there is a trivial forbidden thread $p_{v}$ for some $v \in
Q_0$, we know there exists $\gamma \in Q_1$ with $s(\gamma)=v$ or
$\beta \in Q_1$ with $e(\beta)=v$, in the first case we define
$\sigma (p_v)=\varepsilon (p_v):=-\sigma (\gamma)$, and in the
second one $\sigma(p_v)=\varepsilon(p_v):=-\varepsilon (\beta)$.
\begin{remark}{\normalfont
It is also possible to define threads by using functions $\sigma$ and $\varepsilon$ and definition of strings in \cite{br87}, but then four trivial stings will be needed, two corresponding to the permitted case and two corresponding to the forbidden case. }
\end{remark}

The considered relations are monomials of length two, so we will
indicate them in the quiver by using doted lines, joining each
pair of arrows which form a relation.

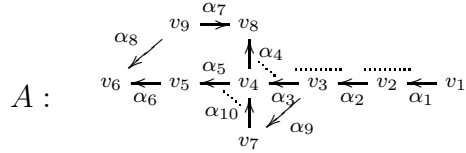
\begin{example}{\normalfont
\[A:\quad
\def\objectstyle{\scriptstyle}
\def\labelstyle{\scriptstyle}
\vcenter{
  \xymatrix@-1.1pc{
      & v_9 \ar[ld]_-{\alpha_8} \ar[r]^-{ \alpha_7}  & v_{8}  \\
      v_6  & v_5 \ar[l]^-{\alpha_6} & v_4 \ar[l]^-{}="d"_-{\alpha_5} \ar[u]^-{}="b"_-{\alpha_4} & v_{3} \ar[l]^-{}="a"_-{\phantom{\cdot}}="f"^-{\alpha_3} \ar[ld]^-{\alpha_9}  & v_2 \ar[l]_-{\phantom{\cdot}}="e"^-{\alpha_2}  & v_{1} \ar[l]_-{\phantom{\cdot}}="g"^-{\alpha_1} \\
           &            & v_7 \ar[u]^-{}="c"^-{\alpha_{10}}       \\
      \ar@{.}"a";"b" \ar@{.} "c";"d" \ar@{.}"e";"f" \ar@{.}"g";"e"
                  }
}
\]

Define $\sigma (\alpha_1)=\sigma (\alpha_4)=\sigma (\alpha_6)=\sigma (\alpha_8)=\sigma (\alpha_9)=\sigma (\alpha_{10})=+1$, $\sigma (\alpha_2)=\sigma (\alpha_3)=\sigma (\alpha_5)=\sigma (\alpha_7)=-1$,  $\varepsilon (\alpha_3)=\varepsilon (\alpha_7)=\varepsilon (\alpha_8)=+1$ and $\varepsilon (\alpha_1)=\varepsilon (\alpha_2)=\varepsilon (\alpha_4)=\varepsilon (\alpha_5)=\varepsilon (\alpha_6)=\varepsilon (\alpha_9)=\varepsilon (\alpha_{10})=-1$.

In this case
$\mathcal{H}_A$ is formed by $ \alpha_1$, $\alpha_4\alpha_{10}\alpha_9\alpha_2$, $\alpha_6\alpha_5\alpha_3$, $\alpha_8$, $\alpha_7$, $1_{v_1}$,$1_{v_7}$ and $1_{v_5}$.

}\end{example}

\subsection{The repetitive algebra of a gentle algebra}\label{repetitiva}

The information of the following Section can be reviewed in
\cite{chj99}, for details see \cite{ri97}. Let $A$ be a
finite-dimensional algebra. The {\em repetitive algebra of $A$} is
defined as the vector space
$$ \hat{A}:=(\oplus_{i\in\mathbb{Z} }A)\oplus(\oplus_{i\in \mathbb{Z}}DA)$$
with multiplication given by
$$(\lambda_i,\phi_i)_{i\in \mathbb{Z}}\cdot(\lambda_i',\phi_i')_{i\in \mathbb{Z}}:=(\lambda_i\lambda_i', \lambda_{i+1}\phi_i'+\phi_i\lambda_i').$$

For a gentle algebra, there is a way to describe its repetitive
algebra constructing a new quiver from an infinite number of
copies of $Q$. More precisely, let $\mathcal{M}$ the set of non
trivial permitted threads of $A$, $M:=\{1,2, \dots,
\#\mathcal{M}\}$ and denote its elements by
$p_i=\alpha_{i,l(p_i)}\dots\alpha_{i,2}\alpha_{i,1}$ for each
$i\in M$. The {\em expansion $\mathbb{Z}(Q, \mathcal{P}  )$} is
defined as the quiver with relations $\mathbb{Z}(Q, \mathcal{P}
):=(\hat{Q},\mathbb{Z}\mathcal{P} )$ where $\hat{Q}$ is given by
$$\hat{Q}_0:=\{v[z]|v\in Q_0, z\in \mathbb{Z} \},$$
$$\hat{Q}_1:=\{\alpha[z]|\alpha\in Q_1, z\in \mathbb{Z} \}\cup\{\alpha_{i,0}[z]|i \in  M, z\in \mathbb{Z}\}$$
with $\alpha[z]:v[z] \rightarrow w[z]$ for each $\alpha:v\rightarrow w,z\in \mathbb{Z}$  and the connection arrows $\alpha_{i,0}[z]:e(\alpha_{i,l(p_i)})[z+1]\rightarrow s(\alpha_{i,1})[z]$ for each $i\in M,z\in\mathbb{Z}$. Let $f:\hat{Q}_1\rightarrow M$ be the function $f(\alpha_{i,j}[z])=i$. Define $$\mathbb{Z}\mathcal{P}:=\{\beta\alpha|\alpha,\beta\in\hat{Q}_1, f(\alpha)\neq f(\beta) \}.$$

There is an automorphism $\nu$ in $(\hat{Q},\mathbb{Z}\mathcal{P} )$ given by $\nu(v[z])=v[z+1]$ and $\nu(\alpha_{i,j}[z])=\alpha_{i,j}[z+1]$,  which corresponds to the Nakayama functor. We say that a path in $\hat{Q}$ not involving any zero relation in $\mathbb{Z}\mathcal{P} $ is {\em full} if it starts in some $v\in \hat{Q}_0$ and ends in $\nu^{-1}(v)$.
Define two sets of relations $$\hat{\mathcal{P}}:=\mathbb{Z}\mathcal{P}\cup \{p-p'|p,p' \text{ are full with $s(p)=s(p')$ and $e(p)=e(p')$}\}$$
$$\phantom{\hat{\mathcal{P}}:=}\cup\{q|q\text{ is a path in $\hat{Q}$ which properly contains a full path}\},$$
$$\bar{\hat{\mathcal{P}}}:=\mathbb{Z}\mathcal{P}\cup \{p|p \text{ is a full path}\}.$$
The special biserial algebra $ {\mathrm k}\hat{Q}/ \left<
\hat{\mathcal{P}} \right>$ is the repetitive algebra of $A$,
$\hat{A}$, while the string algebra $ \bar{\hat{A}}:={\mathrm
k}\hat{Q}/ \left< \bar{\hat{\mathcal{P}}} \right>$ is obtained
from $\hat{A}$ by constructing the quotient over the socles of all
projective-injective modules. For each arrow $\beta$ in
$\bar{\hat{A}}$ there are arrows $\alpha$ and $\gamma$ such that
$\beta\alpha$ and $\gamma\beta$ are permitted paths in $\hat{Q}$
not involving any zero relation in $\bar{\hat{\mathcal{P}}}$, we
say then that $\bar{\hat{A}}$ is {\em expanded}. In this case
there are only two kind of vertices:
\begin{enumerate}
\item
{\em Transition vertices}: There is just one arrow $\alpha$ which ends in the vertex, only one arrow $\beta$ which starts in it and $\beta\alpha\notin \bar{\hat{\mathcal{P}}}$.
\item
{\em Crossing vertices}: There are exactly two arrows which end in the vertex and exactly two arrows which start in it.
\end{enumerate}

As a consequence of the construction of $\bar{\hat{A}}$, a
transition vertex $t$ is the beginning of exactly one non trivial
permitted thread of $\bar{\hat{A}}$, $\mathfrak{p}(t)$, and is the
end of exactly one non trivial permitted thread of
$\bar{\hat{A}}$, $\mathfrak{i}(t)$. Meanwhile, a crossing vertex
is the beginning and ending of exactly two non trivial permitted
threads of $\bar{\hat{A}}$. Also
$$\mathcal{H}_{\bar{\hat{A}}}=\{p|p \text{ is a non trivial permitted thread of $\bar{\hat{A}}$}\}$$
$$\phantom{\mathcal{H}_{\bar{\hat{A}}}}\cup\{t\in\hat{Q}_0|t \text{ is a transition vertex}\}.$$
For each arrow $\beta\in\hat{Q}_1$ define two elements of $\mathcal{H}_{\bar{\hat{A}}}$ in the following way:

$\mathfrak{u}(\beta):=
\begin{cases}
1_{e(\beta)}, \text{ if $e(\beta)$ is a transition vertex}\\
\text{the non trivial permitted thread which ends in $e(\beta)$}\\
\text{ \phantom{rrr}  and does not contain $\beta$, else.}
\end{cases}$

$\mathfrak{v}(\beta):=
\begin{cases}
1_{s(\beta)}, \text{ if $s(\beta)$ is  a transition vertex}\\
\text{the non trivial permitted thread which starts in $s(\beta)$}\\
\text{  \phantom{rrr} and does not contains $\beta$, else.}
\end{cases}$

Therefore, $\mathcal{H}_{\bar{\hat{A}}}$ is the disjoint union of
the sets $\{ \mathfrak{v}(\beta)|\beta\in\hat{Q}_1\}$ and
\\$\{ \mathfrak{p}(t)|t \text{ is  a transition vertex}\}$, or the disjoint union of $\{
\mathfrak{u}(\beta)|\beta\in\hat{Q}_1\}$ and \\ $\{
\mathfrak{i}(t)|t \text{ is a transition vertex}\}$. We can define
a bijection
$\tau:\mathcal{H}_{\bar{\hat{A}}}\rightarrow\mathcal{H}_{\bar{\hat{A}}}$
by
$$\tau(\mathfrak{r}):=
\begin{cases}
\mathfrak{u}(\beta) \text{ if $\mathfrak{r}=\mathfrak{v}(\beta)$}\\
\mathfrak{i}(\nu^{-1}(t)) \text{ if $\mathfrak{r}=\mathfrak{p}(t)$.}
\end{cases}$$

The vertices of the stable Auslander-Reiten quiver of $\hat{A}$
which are end terms of Auslander-Reiten sequences with just one
indecomposable middle term, are naturally parametrized by the set
$\mathcal{H}_{\bar{\hat{A}}}$ and $\tau$ is the action of the
Auslander-Reiten translation, see \cite{br87}. If $Q$ is not a
tree, the infinite $\tau$-orbits of $\mathcal{H}_{\bar{\hat{A}}}$
parametrize the $\mathbb{Z}A_{\infty}$ components of the stable
Auslander-Reiten quiver of $\hat{A}$, while the finite
$\tau$-orbits of $\mathcal{H}_{\bar{\hat{A}}}$ parametrize the
$\mathbb{Z}A_{\infty}/\left< \tau^n \right>$ components which come
from string modules.

\section{Combinatorial calculation of invariants}\label{algoritmo}

Now we describe a combinatorial algorithm to produce certain pairs
of natural numbers, by using only the quiver with relations which
defines a gentle algebra $A$, without needing the Auslander-Reiten
quiver of $\hat{A}$. The description of the modules corresponding
to Auslander-Reiten sequences with just one indecomposable middle
term, see \cite[2.3]{chj99} and \cite[3]{br87}, is the key to
prove that the algorithm provides indeed derived equivalent
invariants. This is presented in the Section \ref{justificacion}.

In the algorithm we go through the quiver which defines the algebra, going forward through permitted threads and backwards through forbidden threads in such a way that each arrow and its inverse is used exactly once.
The procedure is as follows:
\begin{enumerate}
\item
\begin{enumerate}
\item
Begin with a permitted thread of $A$, say $H_0$.
\item
If $H_{i}$ is defined, consider $\Pi_i$, the forbidden thread which ends in $e(H_i)$ and such that $\varepsilon (H_{i})= - \varepsilon (\Pi_i)$.
\item
Let $H_{i+1}$ be the permitted thread which starts in $s(\Pi_i)$ and such that $\sigma (H_{i+1})= - \sigma (\Pi_i)$.
\\
The process stops when $H_n=H_0$ for some natural number $n$. Let $m=\sum_{1\leq i \leq n} l(\Pi_{i-1})$. We obtain the pair $(n,m)$.
\end{enumerate}
\item
Repeat the first step of the algorithm until all permitted threads of $A$ have been considered.
 \item
If there are directed cycles in which each pair of consecutive arrows form a relation, we add a pair $(0,m)$ for each of those cycles, where $m$ is the length of the cycle.
\item
Define $\phi_A :\mathbb{N}^2 \rightarrow \mathbb{N}$ where $\phi_A (n,m)$ is the number of times the pair $(n,m)$ arises in the algorithm.
\end{enumerate}
This function is invariant under derived equivalence, as we shall
see in the Section \ref{justificacion}. Note that $\phi_A$ has
always a finite support. Let $\{(n_1,m_1),(n_2,m_2),\dots
,(n_k,m_k)\}$ be the support of $\phi_A$ , denote $\phi_A$ by
$[(n_1,m_1),(n_2,m_2),\dots ,(n_k,m_k)]$ where each $(n_j,m_j)$ is
written $\phi_A(n_j,m_j)$ times and the order in which they are
written is arbitrary. Define also $\# \phi_A :=\sum _{1\leq j \leq
k } \phi_A(n_j,m_j)$.

\begin{remark}{\normalfont
The case where there are directed cycles in which each pair of
consecutive arrows form a relation, corresponds to algebras of
infinite global dimension. According to the previous algorithm,
for each non trivial forbidden thread $\Pi$ determined by such a
cycle we add a pair $(0,m)$, where $m=l(\Pi)$. This is indeed a
degenerate form of the first step of the process because we go
backwards through this forbidden thread of length $m$ but there is
no going forward through any permitted thread. }
\end{remark}
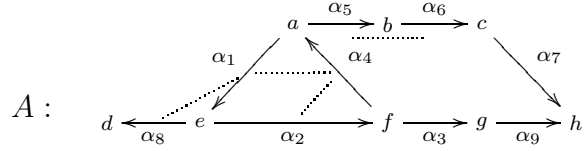
\begin{example}{\normalfont
Consider
\[ A:\quad
\def\objectstyle{\scriptstyle}
\def\labelstyle{\scriptstyle}
\vcenter{
  \xymatrix{
         &                 & a\ar[dl]^-{}="b"_-{\alpha_1}\ar[r]  _-{\phantom{\cdot}}="d" ^-{\alpha_5} & b \ar[r] _-{\phantom{\cdot}}="f"^-{\alpha_6} & c \ar[dr]^-{}="g"^-{\alpha_7}\\
      d  &  e \ar[l]^-{}="c"^-{\alpha_8} \ar[rr]^-{}="e"_-{\alpha_2} &            & f\ar[ul] ^-{}="a"_-{\alpha_4}\ar[r] _-{\alpha_3}& g\ar[r]^-{}="h" _-{\alpha_9} &h  \\
\ar@{.}"a";"b"\ar@{.}"a";"e"\ar@{.}"c";"b"\ar@{.}"d";"f"
                  }
}
\]
Let $\sigma (\alpha_1)=\sigma (\alpha_2)=\sigma (\alpha_3)=\sigma (\alpha_7)=\sigma (\alpha_9)=+1$, $\sigma (\alpha_4)=\sigma (\alpha_5)=\sigma (\alpha_6)=\sigma(\alpha_8)=-1$,  $\varepsilon (\alpha_4)=\varepsilon (\alpha_7)=+1$ and $\varepsilon (\alpha_1)=\varepsilon (\alpha_2)=\varepsilon (\alpha_3)=\varepsilon (\alpha_5)=\varepsilon (\alpha_6)=\varepsilon (\alpha_{8})=\varepsilon (\alpha_9)=-1$.

Let $H_0=\alpha_5\alpha_4$ and $\Pi_0=1_b$ be the trivial forbidden thread in $b$. Then $H_1$ is the permitted thread which starts in $b$ and $\sigma (H_{1})= - \sigma (\Pi_0)$, that is, $\alpha_7\alpha_6$. Now $\Pi_1=\alpha_9$ because it is the forbidden thread which ends in $h$ and $\varepsilon (\Pi_1)=-\varepsilon (H_1)$. Then $H_2$, the permitted thread which starts in $s(\Pi_1)$ and such that $\sigma (H_{2})= - \sigma (\Pi_1)$ is the trivial path in $g$, $1_g$. So, $\Pi_2$, the forbidden thread which ends in $g$ with $\varepsilon (\Pi_2)=-\varepsilon (H_2)$ is $\alpha_3$. Then $H_3=\alpha_5\alpha_4=H_0$, $n=3$ and $m= l(\Pi_0)+l(\Pi_1)+ l(\Pi_2)=0+1+1=2$. The corresponding pair is $(3,2)$. We can write it as follows:
\[
\begin{matrix}
H_0  &= & \alpha_5\alpha_4 & \Pi_0^{-1} & = & 1_b \\
H_1  &= &  \alpha_7\alpha_6  & \Pi_1^{-1} & = &\alpha_9^{-1} \\
H_2  &= & 1_g             & \Pi_2^{-1} & = &\alpha_3^{-1} \\
H_3  &= & H_0 \\
     &  &   & & \rightarrow &(3,2)
\end{matrix}
\]

If we continue with the algorithm we obtain two other pairs,
$(2,4)$ and $(2,3)$, in the following way:
\[
\begin{matrix}
H_0 &= & \alpha_8 & \Pi_0^{-1} & = & 1_d \\
H_1 & = & 1_d       & \Pi_1^{-1} & = &\alpha_2^{-1}\alpha_4^{-1} \alpha_1^{-1}\alpha_8^{-1}  \\
H_2 & = & H_0 \\
 &  &   & & \rightarrow &(2,4) \\
\\
H_0 & = & \alpha_9\alpha_3 \alpha_2\alpha_1 & \Pi_0^{-1} & = &  \alpha_7^{-1} \\
H_1 & = & 1_c                         & \Pi_1^{-1} & = & \alpha_5^{-1}\alpha_6^{-1} \\
H_2 & = & H_0 \\
 &  &   & & \rightarrow &(2,3)\\
\end{matrix}
\]
In this case $\phi_A=[(3,2),(2,4),(2,3)]$.

}\end{example}
\begin{remark}\label{flechashilos}{\normalfont
Let  $\phi_A=[(n_1,m_1),(n_2,m_2),\dots ,(n_k,m_k)]$. As these pairs are the ones obtained by using the algorithm we know that $\sum _{1\leq j \leq k } n_j=\#\mathcal{H}_A$ because each permitted thread appears exactly once. Moreover $\sum _{1\leq j \leq k } m_j=\#Q_1$, because this corresponds to the sum of lengths of all forbidden threads of $A$ and each arrow arises exactly in one of such threads.
In other words $\sum_{(n,m)\in\mathbb{N}^2}\phi_A(n,m)p_1(n,m)=\#\mathcal{H}_A$ and $\sum_{(n,m)\in\mathbb{N}^2}\phi_A(n,m)p_2(n,m)=\#Q_1$, where $p_i:\mathbb{N}^2\to \mathbb{N}$ is the projection in the coordinate $i$, for $i\in \{1,2\}$.
}
\end{remark}
\begin{remark}\label{alg1}{\normalfont
In order to prove that $\phi_A$ is invariant under derived equivalence, let us study all possible cases arising during the algorithm:
\begin{enumerate}
\item
For $H_i$ trivial, $H_i=1_v$ with $v\in Q_0$:
\begin{enumerate}
\item
If $v=s(h)$ for some $h$ non trivial permitted thread of $A$, $\Pi_i$ is trivial and $H_{i+1}=h$.
\item
If $v$ is not the starting point of a non trivial permitted thread of $A$,  $\Pi_i$ is the forbidden thread which ends in $\gamma$, the only arrow in $Q$ such that $e(\gamma)=v$.
\end{enumerate}
Graphically:
\begin{enumerate}
\item
$
\def\objectstyle{\scriptstyle}
\def\labelstyle{\scriptstyle}
\vcenter{
  \xymatrix@-1.1pc{
        & \ar@{.}[l] &        &         & v \ar[lll]^-{}="b"_-{h}
                  }
}
$
\item
$
\def\objectstyle{\scriptstyle}
\def\labelstyle{\scriptstyle}
\vcenter{
  \xymatrix@-1.1pc{
       v  &    \ar[l]^-{}="b"_-{\gamma} &  \ar@{.}[l]
                  }
}
$ or
$
\def\objectstyle{\scriptstyle}
\def\labelstyle{\scriptstyle}
\vcenter{
  \xymatrix@-1.1pc{
   & \ar@{.}[l] & v \ar[l]  &    \ar[l]^-{}="b"_-{\gamma} &  \ar@{.}[l]&.
                  }
}
$
\end{enumerate}
\item
For $H_i$ non trivial:
\begin{enumerate}
\item
If there is $\gamma \in Q_1$ such that $e(H_i)=e(\gamma)$ with $\varepsilon (H_i)=-\varepsilon (\gamma)$, $\Pi_i$ is the only forbidden thread which final arrow is $\gamma$.
\item
If there is no $\gamma \in Q_1$ such that $e(H_i)=e(\gamma)$ with $\varepsilon (H_i)=-\varepsilon (\gamma)$, $\Pi_i$ is trivial and $H_{i+1}$ is the permitted thread such that $s(H_{i+1})=e(H_{i})$, which is $1_{e(H_i)}$ if $e(H_i)$ is a one degree vertex.

\end{enumerate}
This corresponds to:
\begin{enumerate}
\item
$
\def\objectstyle{\scriptstyle}
\def\labelstyle{\scriptstyle}
\vcenter{
  \xymatrix@-1.1pc{
      & & \ar[ld]_-{\gamma} \\
      &\ar@{.}[l] & & & \ar[lll]^-{H_i} & \ar@{.}[l]
                  }
}
$
\item
$
\def\objectstyle{\scriptstyle}
\def\labelstyle{\scriptstyle}
\vcenter{
  \xymatrix@-1.1pc{
   & \ar@{.}[l]  &  & & \ar[lll] ^-{\phantom{\cdot}}="a"_-{H_{i+1}} &    & & \ar[lll]^-{\phantom{\cdot}}="b"_-{H_i} &  \ar@{.}[l]\\
\ar@{.}"a";"b"
                  }
}
$
or
$
\def\objectstyle{\scriptstyle}
\def\labelstyle{\scriptstyle}
\vcenter{
  \xymatrix@-1.1pc{
        & \ar@{}[l]^-{\phantom{\cdot}}="a" &        &         & \ar[lll]^-{\phantom{\cdot}}="b"_-{H_i}&  \ar@{.}[l]&.\\
\ar@{}"a";"b"
                  }
}
$

\end{enumerate}
\end{enumerate}
}\end{remark}

\section{Automorphism groups}\label{automorfismos}

We shall prove that for gentle algebras the number of arrows is
invariant under derived equivalence. This is necessary to define
in a proper way the characteristic component in Section
\ref{invariantes}. From now on we consider ${\mathrm k}$ an
algebraically closed field. By a result presented in \cite{hs01}
and \cite{ro00}, the group $\operatorname{Out}^o(A)$ (which will
be defined later) of a finite-dimensional ${\mathrm k}$-algebra is
invariant under derived equivalence. We proof that for a gentle
algebra it is a group of type $S\ltimes \mathcal{U}$ with $S\cong
({\mathrm k}^*)^{c(Q)}$ and $\mathcal{U}$ a nilpotent subgroup, so
the result follows as a consequence.

For a ${\mathrm k}$-algebra $A$ denote by $\operatorname{Aut}(A)$
its automorphism group, by $\operatorname{Inn}(A)$ its inner
automorphism group, that is $\operatorname{Inn}(A):=\{\iota_x|
x\in A^*\}$ where $\iota _x (a)=x^{-1}ax$ for all $a\in A$, and by
$\operatorname{Out}(A)$ its outer automorphism group,
$\operatorname{Out}(A):= \operatorname{Aut}(A)/
\operatorname{Inn}(A)$. As $\operatorname{Out}(A)$ is an affine
group consider the connected component to which the identity
belongs, called its neutral component and denoted by
$\operatorname{Out}^o(A)$.

\begin{theorem}{\normalfont\label{out}\cite{hs01}, \cite{ro00}}
Let $A$ be a basic finite-dimensional ${\mathrm k}$-algebra. The
affine group $\operatorname{Out}^o(A)$ is invariant under derived
equivalence.
\end{theorem}

In order to make some arguments easier we study a subgroup of $\operatorname{Out}(A)$ which has the same neutral component. Following the article \cite{gs99} define

\[
\begin{aligned}
\operatorname{Aut}^l(A) & :=  \{f\in \operatorname{Aut}(A)|f(v)=v\phantom{*} \forall v \in Q_0\}\\
\operatorname{Inn}^l(A) & :=  \operatorname{Inn}(A)\cap \operatorname{Aut}^l(A) = \{\iota_x\in\operatorname{Inn}(A)|x \in \oplus_{v\in Q_0} vAv \}\\
\operatorname{Out}^l(A) & :=
\operatorname{Aut}^l(A)/\operatorname{Inn}^l(A).
\end{aligned}
\]

We have
\begin{theorem}\label{saorin}{\normalfont\cite [Thm.15]{gs99}}
Let $A$ be a basic finite-dimensional algebra over a field
${\mathrm k}$, then $\operatorname{Out}^l(A)$ is a subgroup of
$\operatorname{Out}(A)$ of finite index and has the same neutral
component $\operatorname{Out}^o(A)$.
\end{theorem}

For a gentle algebra $A={\mathrm k}Q/ \left< \mathcal{P} \right>$
we consider the set of all non zero paths in $A$, that is, which
are not in $\left< \mathcal{P} \right>$,  as a basis. This
collection will be denoted by $\Gamma$; let $\Gamma_{\geq 1}$ be
the set of elements from $\Gamma$ of length greater or equal to
one. Therefore, $A=\oplus_{C\in \Gamma}{\mathrm k}C$ and its
Jacobson radical is $\mathcal{J}=\oplus_{C\in\Gamma_{\geq
1}}{\mathrm k}C$. An element $f$ in $\operatorname{Aut}^l(A)$
fixes the vertices of $Q$, so $f$ is completely determined by the
value it takes in the arrows of the quiver. More precisely, for
$f\in\operatorname{Aut}^l(A)$ and $\alpha\in Q_1$
$$f(\alpha)=\sum_{C\in \Gamma}f_{C}(\alpha)C=\sum_{C\in e(\alpha)\Gamma s(\alpha)}f_{C}(\alpha)C$$
where $f_{C}(\alpha)\in {\mathrm k}$ and the second equality holds
because $f$ fixes the vertices of $Q$. In general, if $D=\alpha_n
\dots\alpha_2\alpha_1$
$$f(D)=\sum_{C_i\in e(\alpha_i)\Gamma s(\alpha_i)}f_{C_n}(\alpha_n)\dots f_{C_2}(\alpha_2)f_{C_1}(\alpha_1)C_n\dots C_2C_1$$
and we denote by $f_C(D)$ the scalar $f_{C_n}(\alpha_n)\dots
f_{C_2}(\alpha_2)f_{C_1}(\alpha_1)$, if $C=C_n\dots C_2C_1$.

Recall that we identify paths with the corresponding classes in
the quotient ${\mathrm k}Q/ \left< \mathcal{P} \right>$; anyway,
the scalars $f_{C}(\alpha)$ are well defined because $A$ is gentle
and consequently monomial, so there are no two different paths
which represent the same element in ${\mathrm k}Q/ \left<
\mathcal{P} \right>$.

Now, we need to study the neutral component of $\operatorname{Out}^l(A)$ in detail. The next result is the key step to do so:

\begin{theorem}\label{orden}
Let $A={\mathrm k}Q/ \left< \mathcal{P} \right>$ be a gentle
algebra, with $Q$ connected not being the Kronecker quiver. There
is a total order $<$ in $\Gamma_{\geq 1}$ such that if
$f\in\operatorname{Aut}^l(A)$ with $f_{\alpha}(\alpha)\neq 0$ for
all $\alpha \in Q_1$,
$$C< D \Rightarrow f_C(D)=0$$
for $C,D\in \Gamma_{\geq 1}$.
\end{theorem}

\proof

We define the total order using the path length, if $C,D\in \Gamma_{\geq 1}$ and $l(C)< l(D) $ then we consider $C<D$. It is moreover possible to order paths of the same length linearly according to the above condition. We proceed in several steps:
\begin{enumerate}
\item We prove that for $C,D\in \Gamma$ with $l(C)< l(D)$,
$f_C(D)=0 $. First consider the case $D=\alpha\in Q_1$,
$C\in\Gamma$ and $\alpha\in Q_1$ with $l(C)<l(\alpha)=1$, that is
$l(C)=0$. As
$$f(\alpha)=\sum_{E\in e(\alpha)\Gamma s(\alpha)}f_{E}(\alpha)E,$$
if $\alpha$ is not a loop, this last expression does not involve $C$ because the paths considered on it start and end in a different vertex. If  $\alpha$ is a loop $e(\alpha)\Gamma s(\alpha)=\{ 1_{s(\alpha)}, \alpha \}$ and $\alpha ^2 =0$ so
\[
\begin{aligned}
0& =f(\alpha ^2)=(f(\alpha ))^2=( f_{1_{s(\alpha)}}(\alpha)1_{s(\alpha)}+f_{\alpha}(\alpha)\alpha ) ^2=\\
 & =(f_{1_{s(\alpha)}}(\alpha))^21_{s(\alpha)}+2f_{1_{s(\alpha)}}(\alpha)f_{\alpha}(\alpha)\alpha
\end{aligned}
\]
Then $(f_{1_{s(\alpha)}}(\alpha))^2=0$ and ${\mathrm k}$ is a
field so this implies $f_{1_{s(\alpha)}}(\alpha)=0$. Anyway $l(C)<
l(\alpha)$ implies $f_C(\alpha)=0$. Consider now
$D=\alpha_n\dots\alpha_2\alpha_1$, then $f(D)=f(\alpha_n)\dots
f(\alpha_2)f(\alpha_1)$ and each $f(\alpha_i)$ is an expression
$$\sum_{C_i\in e(\alpha_i)\Gamma s(\alpha_i)}f_{C_i}(\alpha)C_i.$$

Then
 $$\phantom{***}f(D)=\sum_{C_i\in e(\alpha_i)\Gamma s(\alpha_i)}f_{C_n}(\alpha_n)\dots f_{C_2}(\alpha_2)f_{C_1}(\alpha_1)C_n\dots C_2 C_1.$$
So $f_{C}(D)\neq 0$ if and only if $C=C_n\dots C_2 C_1$ for some $C_i\in e(\alpha_i)\Gamma s(\alpha_i)$ with $f_{C_i}(\alpha_i)\neq 0$ for all $i$; using the previous analysis we know that $l(C_i)\geq l(\alpha)=1$ and then $l(C)\geq n=l(D)$.

We just proved that if $C,D\in \Gamma_{\geq 1}$ with $l(C)< l(D)$ then $f_C(D)=0 $.

Now consider the case of paths of the same length; by the description of $f(C)$ it is enough to study paths $C$ and $D$ such that $s(C)=s(D)$ and $e(C)=e(D)$. First we study what happens for arrows.
\item

Let $\alpha,\beta\in Q_1$ such that $s(\alpha)=s(\beta)$ and $e(\alpha)=e(\beta)$. As $A$ is gentle, connected and not the Kronecker one, suppose without loss of generality there exists $\gamma\in Q_1$ such that $s(\gamma)=e(\alpha)$ with $\gamma\alpha =0$ and $\gamma \beta \neq 0$ (if this is not the case we use analogous arguments for the dual situation). In fact $f_\beta(\alpha)=0$ and then we order $\beta < \alpha$ to fulfill the required condition. Indeed
\[
\begin{aligned}
0 & =f(\gamma\alpha)=f(\gamma)f(\alpha)=\\
  & (\sum_{C\in e(\gamma)\Gamma e(\gamma)}f_{C}(\gamma)C)(\sum_{C\in e(\alpha)\Gamma s(\alpha)}f_{C}(\alpha)C)
\end{aligned}
\]
and the coefficient of $\gamma\beta$ is $f_{\gamma}(\gamma)f_{\beta}(\alpha)$ (by $1.$ no trivial paths are involved in the last expression and the algebra is monomial). By hypothesis $f_{\gamma}(\gamma)\neq 0$, so $f_{\beta}(\alpha)=0$.
\item
Consider now $C,D$ paths of length greater or equal to one, with $s(C)=e(D)$ and $e(C)=e(D)$. In fact $f_{C}(D)=f_{D}(C)=0$ and then we can define the order between them as we like. Let $C=\alpha_n\dots\alpha_2\alpha_1$, and each $f(\alpha_i)$ is an expression $\sum_{E\in e(\alpha_i)\Gamma s(\alpha_i)}f_{E}(\alpha_i)E$, so $f(C)$ is a linear combination of paths $E_n\dots E_2 E_1$ with coefficient $f_{E_n}(\alpha_n)\dots f_{E_2}(\alpha_2) f_{E_1}(\alpha_1)$, each $E_i\in e(\alpha_i)\Gamma s(\alpha_i)$. Let $D=\beta_n\dots\beta_2\beta_1$. If there is an $i\in \{2,\dots ,n\}$ such that $s(\alpha_i)\neq s(\beta_i)$, $D$ does not appear in the linear combination (again using $A$ is monomial). If this is not the case, $s(\alpha_i)= s(\beta_i)$ for all $i\in\{1,2,\dots ,n\}$, by $2.$ $f_{\beta_i}(\alpha_i)=0$ for all $i\in\{1,2,\dots ,n\}$ and then the coefficient of $D$ in the linear combination, $f_{\beta_n}(\alpha_n)\dots f_{\beta_2}(\alpha_2) f_{\beta_1}(\alpha_1)$, is zero, that is, $f_D(C)=0$.
Anyway $f_D(C)=0$, and by symmetry $f_C(D)=0$.

\end{enumerate}

\endproof

Now consider elements $f$ of $\operatorname{Aut}^l(A)$ such that $f_\alpha (\alpha) \neq 0$ for all $\alpha\in Q_1$, that is

$$\operatorname{Aut}_*^l(A):=\{f \in \operatorname{Aut}^l(A) | f_\alpha (\alpha) \neq 0 \phantom{*}\forall \alpha\in Q_1\} .$$

As a consequence of the previous result it is possible to describe
this set:

\begin{corollary}
Let $A={\mathrm k}Q/ \left< \mathcal{P} \right>$ be a gentle
algebra, with $Q$ connected not being the Kronecker quiver.
\begin{enumerate}
\item
The open set $\operatorname{Aut}_*^l(A)$ is  a subgroup of finite index of $\operatorname{Aut}^l(A)$, isomorphic to a lower triangular matrices subgroup.
\item
$\operatorname{Aut}_*^l(A)=S(Q_1)\ltimes N(A)$ where
\begin{equation*}
\begin{array}{llll}
S(Q_1)&:=&\{f \in \operatorname{Aut}_*^l(A) | f (\alpha) \in {\mathrm k}^*\alpha\text{\phantom{ }} \forall \alpha\in Q_1\}\cong ({\mathrm k}^*)^{\#Q_1}&  and \\
N(A)  &:=&\{f \in
\operatorname{Aut}_*^l(A)|f_\alpha(\alpha)=1\text{\phantom{ }}
\forall \alpha\in Q_1\}.&
\end{array}
\end{equation*}
\end{enumerate}
\end{corollary}

\proof
\begin{enumerate}
\item By Theorem \ref{orden} there exists a total order in
$\Gamma_{\geq 1}$ such that the matrices describing the elements
of $\operatorname{Aut}_*^l(A)$ in terms of the corresponding
ordered basis of paths in $A$, are invertible and lower
triangular. Then this is a subgroup of $\operatorname{Aut}^l(A)$,
because the product of two such matrices is of the same type.
\item As $A$ is monomial, each element $(c_\alpha)_{\alpha \in
Q_1}\in({\mathrm k}^*)^{\#Q_1}$ defines an automorphism $f$ in
$\operatorname{Aut}_*^l(A)$ given by $f(\alpha)=c_{\alpha}\alpha$
for each $\alpha\in Q_1$  and this correspondence provides an
isomorphism between $({\mathrm k}^*)^{\#Q_1}$ and $S(Q_1)$. Now
lets see $\operatorname{Aut}_*^l(A)=S(Q_1)\ltimes N(A)$. Let $f\in
\operatorname{Aut}_*^l(A)$. For each $\alpha\in Q_1$,
$f(\alpha)=c_{\alpha}\alpha +\sum_{C\in e(\alpha)\Gamma
s(\alpha)\setminus \{\alpha\}}f_{C}(\alpha)C$. Consider $g\in
S(Q_1)$ defined by $g(\alpha)=c_{\alpha}\alpha$ for all $\alpha\in
Q_1$. But $fg^{-1}\in N(A)$ because
\[
\begin{aligned}
fg^{-1}(\alpha)& =f(g^{-1}(\alpha))=f(c_{\alpha}^{-1}\alpha)=\\
 & =c_{\alpha}^{-1}f(\alpha)=1\alpha +\sum_{C\in e(\alpha)\Gamma s(\alpha)\setminus\{\alpha\}}c_{\alpha}^{-1}f_{C}(\alpha)C
\end{aligned}
\]
then $f=fg^{-1}g\in N(A)S(Q_1)$. Also, if $f\in N(A)\cap S(Q_1)$,
$f(\alpha)=1\alpha$ for all $\alpha\in Q_1$, and then $f=id$. Lets
proof $N(A)\triangleleft\operatorname{Aut}_*^l(A)$. Consider $g\in
S(Q_1)$ such that $g(\alpha)=c_{\alpha}\alpha$ with $c_{\alpha}\in
{\mathrm k}^*$ for all $\alpha\in Q_1$, and $h\in  N(A)$ such that
$$h(\alpha)=1\alpha +\sum_{C\in e(\alpha)\Gamma s(\alpha)\setminus \{\alpha \}}h_{C}(\alpha)C.$$
Then
\[
\begin{aligned}
g^{-1}hg(\alpha) & = g^{-1}h(c_{\alpha}\alpha)=c_{\alpha}g^{-1}h(\alpha)=\\
                 & = c_{\alpha}g^{-1}(1\alpha +\sum_{C\in e(\alpha)\Gamma s(\alpha)\setminus\{\alpha\}}h_{C}(\alpha)C)=\\
                 & = c_{\alpha}(1g^{-1}(\alpha) +\sum_{C\in e(\alpha)\Gamma s(\alpha)\setminus \{\alpha\}}h_{C}(\alpha)g^{-1}(C))=\\
                 & = 1\alpha +\sum_{C\in e(\alpha)\Gamma s(\alpha)\setminus \{\alpha\}}h_{C}(\alpha)c_{\alpha} g^{-1}(C)
\end{aligned}
\]
which is an element of $N(A)$ because if $g^{-1}(C)=\alpha$, $C=g(\alpha)=c_{\alpha}\alpha$ and then $C=\alpha$. So $S(Q_1)$ normalizes $N(A)$. Consider now $f\in \operatorname{Aut}_*^l(A)$ and $h\in N(A)$; $f$ can be written as $f=\tilde{h}g$ for some $\tilde{h}\in N(A)$ and $g\in S(Q_1)$. Then
$$f^{-1}hf=(\tilde{h}g)^{-1}h(\tilde{h}g)=g^{-1}\tilde{h}^{-1}h\tilde{h}g=g^{-1}(\tilde{h}^{-1}h\tilde{h})g$$
which is an element of $N(A)$ because $\tilde{h}^{-1}h\tilde{h}\in N(A)$ and $S(Q_1)$ normalizes $N(A)$. We conclude $N(A)\trianglelefteq\operatorname{Aut}_*^l(A)$.
This leads us to $\operatorname{Aut}_*^l(A)=S(Q_1)\ltimes N(A)$.

\end{enumerate}

\endproof

\begin{theorem}\label{torus}
Let $A={\mathrm k}Q/ \left< \mathcal{P} \right>$ be a gentle
algebra with $Q$ connected  and different of the Kronecker quiver.
$\operatorname{Out}^o(A)$ is isomorphic to a subgroup of lower
triangular invertible matrices which has the form $S\ltimes
\mathcal{U}$ with $S\cong ({\mathrm k}^*)^{c(Q)}$ a maximal torus
and $\mathcal{U}$ a nilpotent subgroup of
$\operatorname{Out}^o(A)$.
\end{theorem}

\proof

Let $x\in \oplus_{v\in Q_o} vAv$ be invertible. Its component in
$v$ is $x_v= c_v 1_v+\sum_{C\in v\Gamma_{\geq 1} v}c_{C}C$ with
$c_v,c_C\in {\mathrm k}$. As the element $x$ is invertible,
$c_v\neq 0$ for all $v\in Q_0$, and then $x_v= c_v 1_v(
1_v+c_v^{-1}\sum_{C\in v\Gamma_{\geq 1} v}c_{C}C)$ with
$c_v^{-1}\sum_{C\in v\Gamma_{\geq 1} v}c_{C}C$ a nilpotent element
which we denote by $n_v\in vAv$ and $c_v\in {\mathrm k}^*$. Define
$s:=\sum _{v\in Q_0} (c_v 1_v)$ and $u:=\sum _{v\in Q_0}
(1_v+n_v)$ observe $su=x$ because $c_v 1_v 1_w=0=c_v 1_v n_w$ for
$v,w\in Q_0$ with $v\neq w$ and then $su =\sum_{v\in Q_0}x_v=x$;
so $\iota _x =\iota_u\iota_s$ and $\iota_u\in N(A)$, $\iota_s\in
S(Q_1)$.

Define
\[
\begin{aligned}
\operatorname{Inn}_N^l(A)& :=\{\iota_u | u=\sum _{v\in Q_0} (1_v+n_v), n_v\in vAv \text{ nilpotent} \}\\
\operatorname{Inn}_S^l(A)& :=\{\iota_s | s=\sum _{v\in Q_0} (c_v
1_v), c_v\in {\mathrm k}^* \}.
\end{aligned}
\]

For $f\in \operatorname{Aut}_*^l(A)$, $\iota_x\in\operatorname{Inn}^l(A)$, $y\in A$
$$(f\iota_x f^{-1})(y)=f(\iota_x f^{-1}(y))=f(x^{-1}f^{-1}(y)x)=f(x^{-1})yf(x)=\iota_{f(x)}(y)$$
and then $f\iota_x f^{-1}=\iota_{f(x)}$. Some remarks can be made at this point:
\begin{enumerate}
\item $\operatorname{Inn}_S^l(A)\trianglelefteq S(Q_1)$ because
$\operatorname{Inn}_S^l(A)\subset Z(\operatorname{Aut}^l(A))$.
Indeed, if $f\operatorname{Aut}^l(A)$ and $\iota_s \in S(Q_1)$
then $f\iota_s f^{-1}=\iota_{f(s)}=\iota_{s}$ because $s$ can be
expressed as $\sum _{v\in Q_0} (c_v 1_v)$, with $ c_v\in {\mathrm
k}^* $ and then $f(s)=s$. \item
$\operatorname{Inn}^l(A)\trianglelefteq \operatorname{Aut}_*^l(A)$
because if $f\in \operatorname{Aut}_*^l(A)$ and $x\in \oplus_{v\in
Q_o} vAv$ $f(x)\in \oplus_{v\in Q_o} vAv$, so $f\iota_x
f^{-1}=\iota_{f(x)}\in\operatorname{Inn}^l(A)$. \item
$\operatorname{Inn}_N^l(A)\trianglelefteq N(A)$ because in fact
$\operatorname{Inn}_N^l(A)$ is a normal subgroup of
$\operatorname{Aut}_*^l(A)$. In order to explain this let $u=\sum
_{v\in Q_0} (1_v+n_v)$ with $n_v= c_v^{-1}\sum_{C\in v\Gamma_{\geq
1} v}c_{C}C$ a nilpotent element and $f\in
\operatorname{Aut}_*^l(A)$; then $f(n_v)=c_v^{-1}\sum_{C\in
v\Gamma_{\geq 1} v}c_{C}f(C)$ which is again nilpotent because no
trivial paths appear in $f(C)$ and $A$ is gentle. Therefore,
$f(u)=\sum _{v\in Q_0} (1_v+\hat{n}_v)$ with $\hat{n}_v$
nilpotent,  and then $f\iota_u f^{-1}=\iota_{f(u)}\in
\operatorname{Inn}_N^l(A)$.
\end{enumerate}

It is possible to consider the morphism
$$\varphi :S(Q_1)\rightarrow \operatorname{Aut}(N(A))$$
for each $g\in S(Q_1)$ $\varphi(g)$ is defined by $\varphi(g)(h):=g^{-1}hg$ for all $h\in N(A)$. This induces a morphism
$$\bar{\varphi} :S(Q_1)/ \operatorname{Inn}_S^l(A)\rightarrow \operatorname{Aut}(N(A)/\operatorname{Inn}_N^l(A)).$$

Also define
$$\psi:S(Q_1)\ltimes N(A)\rightarrow  S(Q_1)/\operatorname{Inn}_S^l(A)\ltimes_{\bar{\varphi}} N(A)/\operatorname{Inn}_N^l(A)$$
by $\psi(gh)=(g\operatorname{Inn}_S^l(A),h\operatorname{Inn}_N^l(A))$. This is a group morphism: if $g,\hat{g}\in S(Q_1)$ and $h,\hat{h}\in N(A)$

\[
\begin{aligned}
\psi((gh)(\hat{g}\hat{h}))& =\psi(g\hat{g}(\hat{g}^{-1}h\hat{g})\hat{h})=(g\hat{g}\operatorname{Inn}_S^l(A),(\hat{g}^{-1}h\hat{g})\hat{h}\operatorname{Inn}_N^l(A))\\
&=(g\hat{g}\operatorname{Inn}_S^l(A),\varphi(g)(h)\hat{h}\operatorname{Inn}_N^l(A))\\
&=(g\operatorname{Inn}_S^l(A),h\operatorname{Inn}_N^l(A))(\hat{g}\operatorname{Inn}_S^l(A),\hat{h}\operatorname{Inn}_N^l(A))\\
&=\psi(gh)\psi(\hat{g}\hat{h})
\end{aligned}
\]
and $\operatorname{ker}\psi=\operatorname{Inn}^l(A)$. This last
statement is because the kernel elements are of the form $gh$ with
$g\in\operatorname{Inn}_S^l(A)$ and $h\in
\operatorname{Inn}_N^l(A)$ (that is $g=\iota_s$ where $s=\sum
_{v\in Q_0} (c_v 1_v)$, $c_v\in {\mathrm k}^*$ and $h=\iota_u$
where $u=\sum _{v\in Q_0} (1_v+n_v)$, $n_v\in vAv$ nilpotent) and,
as shown in the first part of the proof,
$gh=\iota_s\iota_u=\iota_{su}$ with $su\in \oplus_{v\in Q_o} vAv$
an invertible element; in fact any invertible element in
$\oplus_{v\in Q_o} vAv$ can be expressed like this. So

\[
\begin{aligned}
\operatorname{Aut}_*^l(A)/\operatorname{Inn}^l(A)& =(S(Q_1)\ltimes N(A))/\operatorname{Inn}^l(A)\\
 & \cong S(Q_1)/\operatorname{Inn}_S^l(A)\ltimes_{\bar{\varphi}} N(A)/\operatorname{Inn}_N^l(A).
\end{aligned}
\]

As $N(A)$ is nilpotent,  the same is true for
$N(A)/\operatorname{Inn}_N^l(A)$. Consider a vertex $v_0\in Q_0$,
if $s=\sum _{v\in Q_0} (c_v 1_v)$,  $c_v\in {\mathrm k}^*$ and
$s':=1_{v_0}+\sum _{v\in Q_0\setminus \{v_0\}} (c_v' 1_v)$,
$c_v':=c_{v_0}^{-1}c_v$, see that $\iota_s=\iota_{s'}$. Consider
now another element of type $s''=1_{v_0}+\sum _{v\in Q_0\setminus
\{v_0\}} (c_v'' 1_v)$ for some $c_v''\in {\mathrm k}^*$,
$\iota_{s'}=\iota_{s''}$  implies $c_v'=c_v''$ for all $v\in
Q_0\setminus \ \{v_0\}$ because $Q$ is connected. This means that
each of the automorphisms in $\operatorname{Inn}_S^l(A)$ can be
determined by the choice of scalars $c_v$ in the vertices
different from $v_0$ and then   $\operatorname{Inn}_S^l(A)\cong
({\mathrm k}^{*})^{\#Q_0 - 1}$. Therefore,
$S(Q_1)/\operatorname{Inn}_S^l(A)\cong ({\mathrm
k}^{*})^{\#Q_1-(\#Q_0 - 1)}=({\mathrm k}^{*})^{c(Q)}$. By Theorem
\ref{saorin},
$\operatorname{Out}^{o}(A)=(\operatorname{Out}^{l})^{o}(A)$. Also
$\operatorname{Aut}^{l}_*(A)$ is by definition an open set of
$\operatorname{Aut}^{l}(A)$ and using \cite[7.4]{hu95} we know it
is closed at the same time, so
$(\operatorname{Aut}^{l}(A))^{o}=\operatorname{Aut}^{l}_*(A)$.
Then
$(\operatorname{Out}^{l}(A))^{o}=(\operatorname{Aut}^{l}(A))^{o}/\operatorname{Inn}^{l}(A)=\operatorname{Aut}_*^{l}(A)/\operatorname{Inn}^{l}(A)$
and the proof is completed.
\endproof

Now, we are ready to prove Proposition B stated in the
Introduction.
%\begin{corollary}\label{aristasinv}
%The number of cycles and the number of arrows are invariant under derived equivalence for gentle algebras.
%\end{corollary}

By Theorem \ref{out}, $\operatorname{Out}^{o}(A)$ is invariant
under derived equivalence. By Theorem \ref{torus},
$\operatorname{Out}^{o}(A)$ is solvable and contains a maximal
torus of rank $c(Q)$. Finally, for a solvable group all maximal
tori are conjugated and have in particular the same rank. We
conclude that $c(Q)$ is invariant under derived equivalence. As
$\#Q_1= \#Q_0+c(Q)-1$ and the number of vertices is a derived
invariant, so is the number of arrows and this completes the
proof.

\section{Theoretic interpretation of invariants}\label{invariantes}

Let $A={\mathrm k}Q/ \left< \mathcal{P} \right>$ be a gentle
algebra, $Q$ not being a tree and $\hat{A}$ its repetitive
algebra. Recall $\hat{A}$ is selfinjective. Consider the stable
Auslander-Reiten quiver of $\hat{A}$, which we denote by
$\Gamma_{\hat{A},s}$. Let $\tau$ be the Auslander-Reiten
translation and $\Omega^{-1}$ the Heller suspension functor in the
stable module category $\hat{A}$-$\underline{\operatorname{mod}}$;
as $\tau$ and $\Omega$ commute, $\Omega$ permutes the
$\Gamma_{\hat{A},s}$ components.  We call {\em characteristic
components of $\Gamma_{\hat{A},s}$} the components of type
$\mathbb{Z}A_{\infty}$ or $\mathbb{Z}A_{\infty}/\left< \tau^n
\right>$ which come from string modules, see \cite[2.3]{chj99}.
Note that each component $\mathbb{Z}A_{\infty}$ and
$\mathbb{Z}A_{\infty}/\left< \tau^n \right>$, $n\geq 2$ come from
string modules, just the homogeneous tubes
$\mathbb{Z}A_{\infty}/\left< \tau \right>$ consist possibly of
band modules.

We define an equivalence relation on the set of characteristic components of $\Gamma_{\hat{A},s}$ as follows: given $C_1$ and $C_2$ characteristic components of $\Gamma_{\hat{A},s}$, $C_1$ is related to $C_2$ if they are in the same $\Omega$-orbit.
We call an equivalence class under relation a {\em series of components}.

For each series of components $[C]$ there is an associated pair of
natural numbers $(n,m)$, denoted by $i_{[C]}$, obtained as
follows. By  \cite[2.3]{chj99} $\Gamma_{\hat{A},s}$ contains a
finite number of components $\mathbb{Z}A_{\infty}$, and therefore
each of these components $[C]$ has a finite number of elements,
define $i_{[C]}:=(n,m)$ where $|n-m|=\# [C]$ and $\Omega
^{n-m}(M)=\tau ^{n}(M)$ for all module $M$ whose isomorphic class
belongs to this series of component. For a series of components
$[C]$ with $C$ of type $\mathbb{Z}A_{\infty}/\left< \tau^n
\right>$, consider $i_{[C]}:=(n,n)$, this is in fact a limit case
of the previous case. The series consists of a numerable set of
components and $\Omega ^{0}(M)=M=\tau ^{n}(M)$ for all modules $M$
whose isomorphic class belongs to that series, because the
components are tubes or rank $n$. These pairs of natural numbers
and how many times they occur will be essential in this analysis.
Define $N_A:\mathbb{N}^2\rightarrow \mathbb{N}$ with $$N_A(n,m)=\#
\{ [C] \text{ series of components of $\Gamma_{\hat{A},s}$ } |
i_{[C ]}=(n,m)  \}.$$ Just like for the function $\phi_A$
presented in Section \ref{algoritmo}, the support of $N_A$ is
always finite. Let $\{(n_1,m_1),(n_2,m_2),\dots ,(n_k,m_k)\}$ be
the support of $N_A$ , denote $N_A$ by $[(n_1,m_1),(n_2,m_2),\dots
,(n_k,m_k)]$ where each $(n_j,m_j)$ is written $N_A(n_j,m_j)$
times and the order in which they are presented is arbitrary. Also
define $\# N_A :=\sum _{1\leq j \leq k } N_A(n_j,m_j)$ or the
number of series of components of $A$.

\section{Justification of the algorithm}\label{justificacion}

According to Section \ref{repetitiva}, for a gentle algebra, a
module in the stable Auslander-Reiten quiver of the repetitive
algebra which is the start or end of an Auslander-Reiten sequence
with just one indecomposable middle term, corresponds to a
permitted thread of $\bar{\hat{A}}$;         identifying this
modules with the related permitted threads, now we analyze the
action of the Heller suspension functor over the set of permitted
threads. By definition, if $M$ is a module over an algebra $A$,
$\Omega^{-1}(M)=\operatorname{coker}(M\rightarrow I )$ where $I$
is the injective envelope of $M$. Let $H$ be a permitted thread of
$\bar{\hat{A}}$. $H$ is trivial or is obtained from a full path by
deleting its final arrow, which according to the notation of
Section \ref{repetitiva}, can be written like $\alpha_{i,j}[z]$
for some $i\in M$, $j\in \{0,1,\ldots , l(H)\}$, $z\in
\mathbb{Z}$, it also is denoted by $\alpha_{i,j}^{H}[z]$ . In this
case, if $j\neq 0$, $\alpha_{i,j}\in Q_1$ and $H$ is identify with
$\alpha_{i,j}^{H}[z]^{-1}$; else, it is a connection arrow  (see
Section \ref{repetitiva}) and then $H$ corresponds to one of the
copies of the permitted thread of $A$, $h$, in the quiver
associated to $\bar{\hat{A}}$, this is,
$H=\alpha_{i,l(H)}[z]\cdots\alpha_{i,2}[z]\alpha_{i,1}[z]$ with
$h=\alpha_{i,l(H)}\cdots\alpha_{i,2}\alpha_{i,1}$ a permitted
thread of $A$ and $z\in \mathbb{Z}$, denote then $H=h[z]$. Because
of the structure of $\bar{\hat{A}}$ we get

$\Omega^{-1}(H)=
\begin{cases}
\mathfrak{p}(v[z+1]) \text{ if $H=1_{v[z]}$ with $v\in Q_0$, $z\in \mathbb{Z}$ }\\
\mathfrak{v}(\alpha_{i,j}^{H}[z+1]) \text{ if $H$ is non trivial}
\end{cases}$

more precisely,

$\Omega^{-1}(H)=
\begin{cases}
\mathfrak{p}(v[z+1]) \text{ if $H=1_{v[z]}$ with $v\in Q_0$, $z\in \mathbb{Z}$ }\\
1_{\nu(e(H))} \text{ if $H$ is non trivial and $\nu(e(H))$ is a transition}\\
\text{\phantom{$1_{\nu(e(H))}$}vertex}\\
\text{the permitted thread starting in $\nu(e(H))$ which does not} \\
\text{\phantom{$1_{\nu(e(H))}$}involves $\alpha_{i,j}^{H}[z+1]$ if $H$ is non trivial and}\\
\text{\phantom{$1_{\nu(e(H))}$}$\nu(e(H))$ is a crossing vertex}
\end{cases}$

where $\mathfrak{p}$ and $\mathfrak{v}$ are defined as in Section
\ref{repetitiva}.

\begin{remark}\label{alg2}{\normalfont
By the previous analysis and the fact that $A$ is gentle $\Omega^{-1}(H)$ can be described for a permitted thread of $\bar{\hat{A}}$, $H$, in all cases:
\begin{enumerate}
\item
For $H$ trivial, $H=1_{v[z]}$ with $v\in Q_0$ , $z\in \mathbb{Z}$ (in this case $1_v$ is a trivial thread of $A$):
\begin{enumerate}
\item
If $v=s(h)$ for some permitted thread of $A$ $h$, $\Omega^{-1}(H)=h[z+1]$.
\item
If $v$ is not the start point of a permitted non trivial thread in $A$,  $\Omega^{-1}(H)$ is identified with $\gamma[z]$ where $\gamma$ is the only arrow in  $Q$ such that $e(\gamma)=v$.
\end{enumerate}
Graphically:
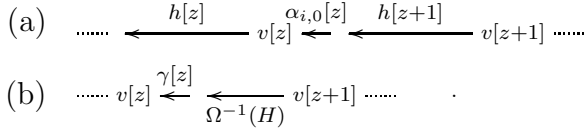
\begin{enumerate}
\item
$
\def\objectstyle{\scriptstyle}
\def\labelstyle{\scriptstyle}
\vcenter{
  \xymatrix@-1.1pc{
        & \ar@{.}[l] &        &         & v[z] \ar[lll]^-{}="b"_-{h[z]} &  \ar[l]^-{}="b"_-{\alpha_{i,0}[z]} & & &v[z+1]\ar[lll]^-{}="b"_-{h[z+1]}& \ar@{.}[l]
                  }
}
$
\item
$
\def\objectstyle{\scriptstyle}
\def\labelstyle{\scriptstyle}
\vcenter{
  \xymatrix@-1.1pc{
      &  v[z]\ar@{.}[l]  &    \ar[l]^-{}="b"_-{\gamma[z]} & &  v[z+1]\ar[ll]^-{ \Omega^{-1}(H)} &  \ar@{.}[l]&.
                  }
}
$
\end{enumerate}
\item
In the case $H=h[z]$ for some permitted non trivial thread of $A$, $h$:
\begin{enumerate}
\item
If there exists $\gamma \in Q_1$ such that $e(h)=e(\gamma)$ with $\varepsilon (h)=-\varepsilon (\gamma)$, $\Omega^{-1}(H)$ is identified with $\gamma[z]^{-1}$, and $\gamma$ corresponds to the final arrow of a forbidden thread of $A$.
\item
If there is no $\gamma \in Q_1$ such that $e(h)=e(\gamma)$ with $\varepsilon (h)=-\varepsilon (\gamma)$, $\Omega^{-1}(H)=h'[z+1]$ where $h'$ is the permitted thread such that $s(h')=e(h)$, which is $1_{\nu(e(H))}$ in case $e(h)$ has degree one.
\end{enumerate}
Graphically:
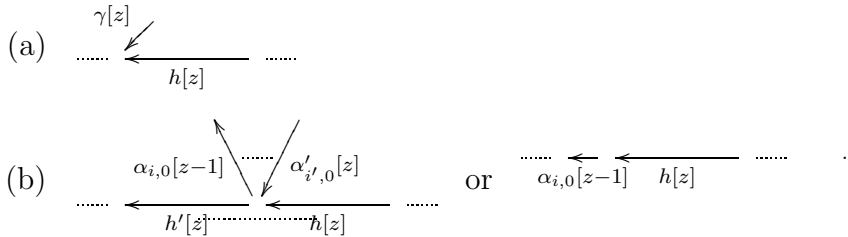
\begin{enumerate}
\item
$
\def\objectstyle{\scriptstyle}
\def\labelstyle{\scriptstyle}
\vcenter{
  \xymatrix@-1.1pc{
      & & \ar[ld]_-{\gamma[z]} \\
      &\ar@{.}[l] & & & \ar[lll]^-{h[z]} & \ar@{.}[l]
                  }
}
$
\item
$
\def\objectstyle{\scriptstyle}
\def\labelstyle{\scriptstyle}
\vcenter{
  \xymatrix@-1.1pc{
    & & &             &                                        &  \ar[ldd]^-{}="c"^-{\alpha_{i',0}'[z]}       \\
    & & &             &                                        &   \\
   & \ar@{.}[l]  &  & & \ar[lll] ^-{\phantom{\cdot}}="a"^-{h'[z]} \ar[uul]^-{}="d"^-{\alpha_{i,0}[z-1]} &    & & \ar[lll]^-{\phantom{\cdot}}="b"^-{h[z]} &  \ar@{.}[l]\\
\ar@{.}"a";"b"\ar@{.}"c";"d"
                  }
}
$
or
$
\def\objectstyle{\scriptstyle}
\def\labelstyle{\scriptstyle}
\vcenter{
  \xymatrix@-1.1pc{
     & \ar@{.}[l] &  \ar[l]^-{\alpha_{i,0}[z-1]}^-{\phantom{\cdot}}="a" &        &         & \ar[lll]^-{\phantom{\cdot}}="b"^-{h[z]}&  \ar@{.}[l]&.\\
\ar@{}"a";"b"
                  }
}
$

\end{enumerate}
\item
In case $H$ is identified with $\alpha_{i,j}^{H}[z]^{-1}$:
\begin{enumerate}
\item If there exists $\gamma \in Q_1$ such that $e(\gamma
[z])=s({\alpha_{i,j}^{H}[z]})=e(H)$ with
$\varepsilon(\gamma)=\sigma (\alpha_{i,j}^{H})$, $\Omega^{-1}(H)$
is identified with $\gamma[z]^{-1}$ and ${\alpha_{i,j}^{H}}\gamma$
is a forbidden thread in $A$. \item If there is no $\gamma \in
Q_1$ such that $e(\gamma[z])=s({\alpha_{i,j}^{H}}[z])=e(H)$ with
$\varepsilon(\gamma)=\sigma (\alpha_{i,j}^{H})$,
$\alpha_{i,j}^{H}$ is the starting point of a forbidden thread in
$A$ and $\Omega^{-1}(H)=h'[z+1]$ where $h'$ is the permitted
thread of $A$ such that $s(h')=s(\alpha_{i,j}^{H})$ and does not
involve $\alpha_{i,j}^{H}$, that is, such that
$\sigma(h')=-\sigma(\alpha_{i,j}^{H})$, which is $1_{\nu(e(H))}$
if $e(H)$ is a transition vertex.
\end{enumerate}
Graphically:
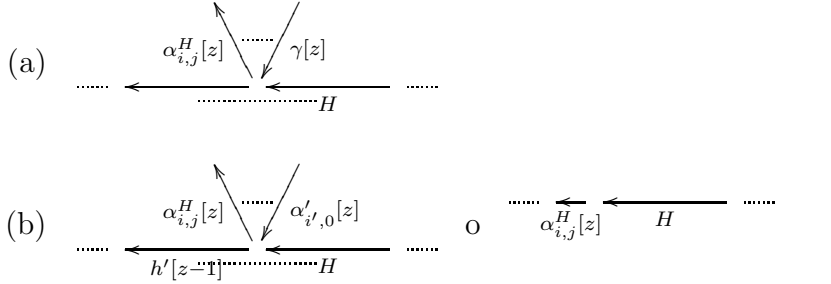
\begin{enumerate}
\item
$
\def\objectstyle{\scriptstyle}
\def\labelstyle{\scriptstyle}
\vcenter{
  \xymatrix@-1.1pc{
    & & &             &                                        &  \ar[ldd]^-{}="c"^-{\gamma[z]}       \\
    & & &             &                                        &   \\
   & \ar@{.}[l]  &  & & \ar[lll] ^-{\phantom{\cdot}}="a"^-{} \ar[uul]^-{}="d"^-{\alpha_{i,j}^{H}[z]} &    & & \ar[lll]^-{\phantom{\cdot}}="b"^-{H} &  \ar@{.}[l]\\
\ar@{.}"a";"b"\ar@{.}"c";"d"
                  }
}
$
\item
$
\def\objectstyle{\scriptstyle}
\def\labelstyle{\scriptstyle}
\vcenter{
  \xymatrix@-1.1pc{
    & & &             &                                        &  \ar[ldd]^-{}="c"^-{\alpha_{i',0}'[z]}       \\
    & & &             &                                        &   \\
   & \ar@{.}[l]  &  & & \ar[lll] ^-{\phantom{\cdot}}="a"^-{h'[z-1]} \ar[uul]^-{}="d"^-{\alpha_{i,j}^{H}[z]} &    & & \ar[lll]^-{\phantom{\cdot}}="b"^-{H} &  \ar@{.}[l]\\
\ar@{.}"a";"b"\ar@{.}"c";"d"
                  }
}
$
o
$
\def\objectstyle{\scriptstyle}
\def\labelstyle{\scriptstyle}
\vcenter{
  \xymatrix@-1.1pc{
     & \ar@{.}[l] &  \ar[l]^-{\alpha_{i,j}^{H}[z]}^-{\phantom{\cdot}}="a" &        &         & \ar[lll]^-{\phantom{\cdot}}="b"^-{H}&  \ar@{.}[l]&.\\
\ar@{}"a";"b"
                  }
}
$
\end{enumerate}
\end{enumerate}
}\end{remark}

After this analysis we get the following result:

\begin{theorem}\label{phi=N}
Let $A$ be a gentle algebra and $\phi_A :\mathbb{N}^2 \rightarrow
\mathbb{N}$ is the function defined in Section \ref{algoritmo},
then $\phi_A=N_A$.
\end{theorem}

\proof

Recall the algorithm of Section \ref{algoritmo}. Comparing Remarks
\ref{alg1} and \ref{alg2} we conclude that, if $H$ is a permitted
thread of $\bar{\hat{A}}$, the powers $\Omega^{i}(H)$ are the
different permitted threads of $A$ appearing in the algorithm, or
identify with the inverses of the arrows in $Q$ which constitute
the non trivial forbidden threads of $A$, obtained in the same
order as in the algorithm. More precisely, let $H_0$ be a
permitted thread of $A$ and $H_0,H_1,\ldots H_n$ the permitted
threads obtained in the process starting with $H_0$. Denote by
$\Pi_0,\Pi_1,\ldots \Pi_{n-1}$ the non trivial forbidden threads
of $A$ involved in this part of the algorithm, where
$$\Pi_i=\pi_{i,l(\Pi_i)}\cdots \pi_{i,2}\pi_{i,1}$$
 with $\pi_{i,j}\in Q_1$ for $i\in\{0,1,2,\ldots,n-1\}$, $j\in\{1,2,\ldots,l(\Pi_i)\}$, $m=\sum_{0\leq i \leq n-1}l(\Pi_i)$ and $(n,m)$ is the pair of natural numbers obtained from $H_0$. This part of the algorithm can be described by the array

\begin{equation*}
\begin{array}{ll}
H_0[0]  &  \pi_{1,1}^{-1}[0]\cdots\pi_{1,l(\Pi_0)}^{-1}[0] \\
H_1[1]  &  \pi_{2,1}^{-1}[1]\cdots \pi_{2,l(\Pi_1)}^{-1}[1] \\
H_2[2]  &  \dots \\
\vdots  &  \\
H_n[n]  &
\end{array}
\end{equation*}
where none of the trivial forbidden threads is written out (and
this is indicated by putting an asterisk instead). Identifying
permitted threads of $\bar{\hat{A}}$ which contain a connection
arrow with its corresponding one $\alpha_{i,j}^{H}[z]^{-1}$ and
using the previous remarks, this is the array
 \begin{equation*}
\begin{array}{ll}
H_0[0]                &  \Omega^{-l(\Pi_0)}(H_0[0])\cdots\Omega^{-1}(H_0[0]) \\
\Omega^{-l(\Pi_0)-1}(H_0[0])  &  \Omega^{-l(\Pi_0)-l(\Pi_1)-1}(H_0[0])\cdots\Omega^{-l(\Pi_0)-2}(H_0[0]) \\
 \Omega^{-l(\Pi_0)-l(\Pi_1)-2}(H_0[0])  &  \dots \\
\vdots  &  \\
 \Omega^{-l(\Pi_0)-l(\Pi_1)-\dots -l(\Pi_{n-1})-n}(H_0[0])  &   \\
\end{array}
\end{equation*}
where the last element is $\Omega^{-m-n}(H_0[0])$ because of the way $m$ was defined. As the algorithm stops in the first moment in which $H_n=H_0$, $\Omega^{-n-m}(H_0[0])$ is the first module in the same $\nu$-orbit of $H_0[0]$ obtained by applying certain number of times $\Omega^{-1}$, in fact $\Omega^{-n-m}(H_0[0])=H_n[n]=\nu^{n}(H_0[0])$. As $\Omega$ is an equivalence, the information about the modules of the Auslander-Reiten quiver which are extreme points of an Auslander-Reiten sequence with one indecomposable middle term, is enough to conclude that $\Omega^{-n-m}=\nu^{n}$ for every module in the series of components where $H_0[0]$ lies. As $\tau=\Omega^{2}\circ\nu$ and $\Omega$ commute with $\nu$ then
$$\tau^{n} = (\Omega^{2}\circ\nu)^{n}   = \Omega^{2n}\circ\nu^{n}\Rightarrow \Omega^{-2n}\circ\tau^{n}  =  \nu^{n}=\Omega^{-n-m}$$
which proves that $\tau^{n}= \Omega^{n-m}$. This is also true if
naturals $n$ and $m$ are equal and in this case $\tau^{n}=
\Omega^{0}=id$ for modules in the series of components where
$H_0[0]$ lies, and $n$ is the least natural for which this
happens, then this components are tubes of rank $n$. If $n\neq m$,
evaluating powers of $\Omega^{-1}$ at $H_0[0]$, we get modules in
the different components forming the series associated to $H_0[0]$
until we get to the same component where we started, after exactly
$n-m$ steps. Then $|n-m|$ is the number of components of the
series. In the case of a gentle algebra of infinite global
dimension, there exists a directed cycle in which every pair of
consecutive arrows define a relation. Let $\alpha$ be any of these
arrows and $H$ the permitted thread in $\bar{\hat{A}}$ with which
we identify $\alpha[0]^{-1}$. By previous analysis, the powers of
$\Omega^{-1}(H), \Omega^{-2}(H)\dots \Omega^{-m}(H)$ identify
precisely with the arrows forming the directed cycle just
mentioned, with $m$ its length and $\Omega^{-m}(H)=H$. This
implies that $\Omega^{-m}=id=\tau^{0}$ for modules in the series
of components where $H$ lies, formed by exactly $m$ components.
\endproof

\begin{example}{\normalfont
For the algebra presented in the example of Section
\ref{algoritmo}, if $H_0 = \alpha_5\alpha_4$ the algorithm is
codified by
  \begin{equation*}
\begin{array}{ll}
(\alpha_5\alpha_4)[0] & *\\
(\alpha_7\alpha_6)[1] &  \alpha_9[1]^{-1} \\
 1_g[2]               &  \alpha_3[2]^{-1} \\
(\alpha_5\alpha_4)[3] &  *
\end{array}
\end{equation*}
for $H_0=\alpha_8 $ by
 \begin{equation*}
\begin{array}{ll}
\alpha_8[0] & *\\
1_d[1]      &  \alpha_2[1]^{-1}\alpha_4[1]^{-1}\alpha_1[1]^{-1}\alpha_8[1]^{-1} \\
\alpha_8[2] &
\end{array}
\end{equation*}
and for $H_0 =\alpha_9\alpha_3\alpha_2\alpha_1$
 \begin{equation*}
\begin{array}{ll}
(\alpha_9\alpha_3\alpha_2\alpha_1)[0] & \alpha_7[0]^{-1}\\
1 _c[1]      &  \alpha_5[1]^{-1}\alpha_6[1]^{-1} \\
(\alpha_9\alpha_3\alpha_2\alpha_1)[2] &
\end{array}
\end{equation*}
which correspond to the following arrays
  \begin{equation*}
\begin{array}{ll}
(\alpha_5\alpha_4)[0] & *\\
\Omega^{-1}((\alpha_5\alpha_4)[0]) &  \Omega^{-2}((\alpha_5\alpha_4)[0]) \\
\Omega^{-3}((\alpha_5\alpha_4)[0]) &  \Omega^{-4}((\alpha_5\alpha_4)[0]) \\
\Omega^{-5}((\alpha_5\alpha_4)[0]) &
\end{array}
\end{equation*}

  \begin{equation*}
\begin{array}{ll}
\alpha_8[0] & *\\
\Omega^{-1}(\alpha_8[0]) &  \Omega^{-5}(\alpha_8[0])\Omega^{-4}(\alpha_8[0])\Omega^{-3}(\alpha_8[0])\Omega^{-2}(\alpha_8[0]) \\
\Omega^{-6}(\alpha_8[0]) &
\end{array}
\end{equation*}
and
  \begin{equation*}
\begin{array}{ll}
(\alpha_9\alpha_3\alpha_2\alpha_1)[0] & \Omega^{-1}((\alpha_9\alpha_3\alpha_2\alpha_1)[0])\\
\Omega^{-2}((\alpha_9\alpha_3\alpha_2\alpha_1)[0]) &  \Omega^{-4}((\alpha_9\alpha_3\alpha_2\alpha_1)[0])\Omega^{-3}((\alpha_9\alpha_3\alpha_2\alpha_1)[0]) \\
\Omega^{-5}((\alpha_9\alpha_3\alpha_2\alpha_1)[0]) & .
\end{array}
\end{equation*}
}\end{example}

\subsection{Main theorem}

We prove now Theorem A stated in the Introduction.
%\begin{theorem}\label{soninv} Let $A$ and
%$B$ be gentle algebras. If $A$ and $B$ are derived equivalent then
%$\phi_A=\phi_B$.
%\end{theorem}
\\
Consider $A={\mathrm k}Q/ \left< \mathcal{P} \right>$ and
$B={\mathrm k}Q'/ \left< \mathcal{P}' \right>$ two derived
equivalent gentle algebras and $\phi_A$ and $\phi_B$ the functions
defined by the pairs of natural numbers obtained in the algorithm
defined in Section \ref{algoritmo}.

As $A$ and $B$ are derived equivalent then $\hat{A}$ and $\hat{B}$
also are, see \cite[Thm. 1.5]{as97} and \cite{as98}, and because
these are selfinjective, by \cite{ric89} we know
$\hat{A}$-$\underline{\operatorname{mod}}$ and
$\hat{A}$-$\underline{\operatorname{mod}}$ are equivalent as
triangulated categories. Then, a series of components
$\mathbb{Z}A_{\infty}$ which has $(n,m)$ as associated pair
corresponds by the equivalence to a series of components of the
same type and equal number of elements, $|n-m|$, and in such a way
that $\Omega ^{n-m}(M)=\tau ^{n}(M)$ for all $M$ a module whose
isomorphism class lies in that series. A series of components of
tubes of rank $n$, with $n\gneq 1$ which has $(n,n)$ as associated
pair, corresponds by the equivalence to a series of components of
tubes of the same rank, to which we associate also the pair
$(n,n)$.

Denote by $\phi_A'$ the restriction of $\phi_A$ to $\mathbb{N}^2\setminus\{(1,1)\}$, as $\phi_A'=N_A'$ by Theorem \ref{phi=N}, $\phi_A'$ describes the action of $\Omega^{-1}$ over the components $\mathbb{Z}A_{\infty}$ and $\mathbb{Z}A_{\infty}/ \left< \tau^n \right>$ of $\hat{A}$ with $n\in\mathbb{N}\setminus\{1\}$ and we conclude then that $\phi_A'=\phi_B'$.

Using Remark \ref{flechashilos} we know that
$$\sum_{(n,m)\in\mathbb{N}^2}\phi_A(n,m)p_2(n,m)=\#Q_1,$$
so $\phi_A(1,1)=\#Q_1-\sum_{(n,m)\in\mathbb{N}^2\setminus
\{(1,1)\}}\phi_A'(n,m)p_2(n,m)$ and similar for $B$. By the last
Section the number of arrows is a derived equivalent invariant so
$\#Q_1=\#Q_1'$, also $\phi_A'(n,m)=\phi_B'(n,m)$, then
$\phi_A(1,1)=\phi_B(1,1)$. Therefore $\phi_A=\phi_B$ and the proof
concludes.

\section{Discussion of the results}\label{clas}

The first aim of this Section is to prove Theorem C, which states
that our invariants separate the derived equivalence classes of
gentle algebras with at most one cycle.
%\begin{theorem}\label{unciclo}
%Let  $A={\mathrm k}Q/ \left< \mathcal{P} \right>$ and $B={\mathrm
%k}Q'/ \left< \mathcal{P}' \right>$ be gentle algebras such that
%$c(Q),c(Q')\leq 1$. Then $A$ and $B$ are derived equivalent if and
%only if $\phi_A=\phi_B$.
%\end{theorem}
\\
Recall, that the derived classification of gentle algebras with at
most one cycle is known \cite{ah81}, \cite{as87}, \cite{dv01} and
\cite{bg04}. More precisely, we have the following three cases:
\begin{enumerate}
\item $Q$ is a tree if and only if $A={\mathrm k}Q/ \left<
\mathcal{P} \right>$ is derived equivalent to $\mathbb{A}_n$ with
$n=\#Q_0$, and so $\phi_A=\phi_{\mathbb{A}_n}=[(n+1,n-1)]$. \item
$Q$ has just one cycle and satisfies the clock condition (that is
the number of clockwise relations equals the number of
anticlockwise relations in the only cycle of $A$) if and only if
$A={\mathrm k}Q/ \left< \mathcal{P} \right>$ is derived equivalent
to a hereditary algebra of type $\tilde{\mathbb{A}}_{p,q}$ for
some $p,q \in \mathbb{N}$ with $p+q=n$. In this case
$\phi_A=\phi_{\tilde{\mathbb{A}}_{p,q}}=[(p,p),(q,q)]$. \item
  $Q$ has just one cycle and $r:=|c-a|\neq 0$ where $c$ and $a$ are the number of clockwise and anticlockwise relations resp. in the only cycle of $A$ if and only if $A={\mathrm k}Q/ \left< \mathcal{P} \right>$ is derived equivalent to an algebra of type $\Lambda (r,n,m)={\mathrm k}Q(r,n,m)/ \left< \mathcal{P}(r,n,m) \right>$ for some $r,n,m \in \mathbb{N}$ with $n\geq r\geq 1$ and $m\geq 0$, graphically:
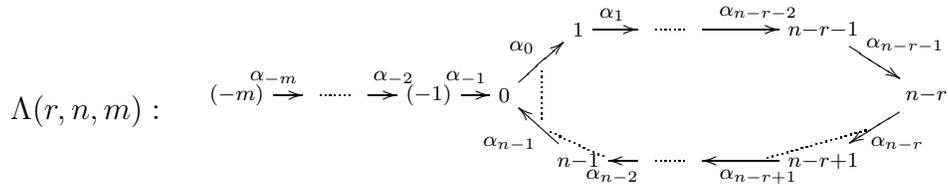
\[\Lambda (r,n,m):\quad
\def\objectstyle{\scriptstyle}
\def\labelstyle{\scriptstyle}
\vcenter{
  \xymatrix@-1.1pc{
  &&&&& 1 \ar[r]^{\alpha_1}&  \ar@{.}[r]  & \ar[rr]^{\alpha_{n-r-2}}& & n-r-1\ar[rd]^{\alpha_{n-r-1}}\\
(-m)\ar[r]^{\alpha_{-m}}&  \ar@{.}[r]  & \ar[r]^{\alpha_{-2}}& (-1)\ar[r]^{\alpha_{-1}} & 0\ar[ru]^{\alpha_{0}}_{}="a"&&&&&& n-r\ar[ld]^{\alpha_{n-r}}_{}="e"\\
 &&&&& n-1 \ar[lu]^{\alpha_{n-1}}_{}="b"&  \ar[l]^{\alpha_{n-2}}_{}="c"  & \ar@{.}[l]& &n-r+1\ar[ll]^{\alpha_{n-r+1}}_{}="d"\\
     \ar@{.} "a";"b" \ar@{.}"b";"c"  \ar@{.}"d";"e"
                  }
}
\]
with $\mathcal{P}(r,n,m)=\conj{ \alpha_0\alpha_{n-1}, \alpha_{n-1}\alpha_{n-2}, \dots , \alpha_{n-r+1}\alpha_{n-r}}$, see \cite{bg04}.
In this case $\phi_A=\phi_{\Lambda (r,n,m)}=[(r+m,m),(n-r,n)]$.
\end{enumerate}
Then if $A={\mathrm k}Q/ \left< \mathcal{P} \right>$ and
$B={\mathrm k}Q'/ \left< \mathcal{P}' \right>$ are gentle algebras
with $c(Q),c(Q')\leq 1$ such that $\phi_A=\phi_B$, they are
derived equivalent. This concludes the proof.

\begin{remark}

By Theorem C our invariants can distinguish two gentle algebras
which are not derived equivalent in all known cases. However,
consider the following two gentle algebras:

$$A:\quad
\def\objectstyle{\scriptstyle}
\def\labelstyle{\scriptstyle}
\vcenter
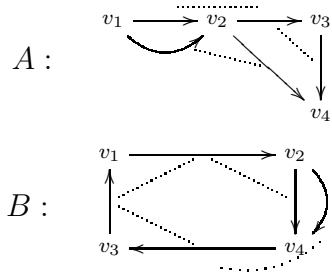
{
  \xymatrix{
     v_1 \ar[r]^{\phantom{\cdot}}="c"\ar@/_1pc/[r]^(.7){}="e"  & v_2 \ar[r]^{\phantom{\cdot}}="d"_{}="g" \ar[dr]_{}="f" & v_3 \ar[d]_{}="h"    \\
                                                              &                                               & v_4
     \ar@{.} "c";"d"\ar@{.}"e";"f" \ar@{.}"g";"h"
                  }
}
$$
$$
B:\quad
\def\objectstyle{\scriptstyle}
\def\labelstyle{\scriptstyle}
\vcenter{
  \xymatrix{
     v_1 \ar[rr]_{}="e"&  &v_2 \ar[d]_{}="f" \ar@/^1pc/[d]^(.7){\phantom{\cdot}}="g" \\
     v_3 \ar[u]_{}="d"&   &v_4 \ar[ll]_{}="c"^{\phantom{\cdot}}="h"
     \ar@{.} "c";"d" \ar@{.}"d";"e" \ar@{.}"e";"f"  \ar@/^.7pc/@{.}"g";"h"
                  }
}
$$

We find $\phi_A=\phi_B=[(3,5)]$. On the other hand we checked with
a computer that it is not possible to transform $A$ into $B$ by a
sequence of elementary derived equivalences from \cite{jan01}. It
should be interesting to find another way to decide if these two
algebras are derived equivalent or not.

\end{remark}

\end{document}